\documentclass[10pt]{article}
\usepackage{amsmath,amssymb,amsfonts}

\newtheorem{theorem}{Theorem}
\newtheorem{lemma}{Lemma}
\renewcommand{\theequation}{\thesection.\arabic{equation}}

\date{}
\title{The Size of the Largest Part of Random Weighted Partitions of Large Integers}
\bigskip

\author{{\bf Ljuben Mutafchiev}\\
American University in Bulgaria, 2700 Blagoevgrad, Bulgaria \\ and
Institute of Mathematics and Informatics of the \\ Bulgarian
Academy of Sciences
\\ \tt {ljuben@aubg.bg}}

\begin{document}
\maketitle

\begin{abstract}
 We consider partitions of the positive integer $n$ whose parts
 satisfy the following condition: for a given sequence of
 non-negative numbers $\{b_k\}_{k\ge 1}$, a part of size $k$
 appears in exactly $b_k$ possible types. Assuming that a weighted
 partition is selected uniformly at random from the set of all
 such partitions, we study the asymptotic behavior of the largest
 part $X_n$. Let $D(s)=\sum_{k=1}^\infty b_k k^{-s}, s=\sigma+iy,$
 be the Dirichlet generating series of the weights $b_k$. Under
 certain fairly general assumptions, Meinardus (1954) has obtained
 the asymptotic of the total number of such partitions as
 $n\to\infty$. Using the Meinardus scheme of conditions, we prove that
 $X_n$, appropriately normalized, converges weakly to a random
 variable having Gumbel distribution (i.e. its distribution
 function equals $e^{-e^{-t}}, -\infty<t<\infty$). This limit
 theorem extends some known results on particular types of
 partitions and on the Bose-Einstein model of ideal gas.

\end{abstract}

\vspace{.5cm}


\vspace{.2cm}

\section{Introduction and Statement of the Result}

A weighted partition of the positive integer $n$ is a multiset of
size $n$ whose decomposition into a union of disjoint components
(parts) satisfies the following condition: for a given sequence of
non-negative numbers $\{b_k\}_{k\ge 1}$, a part of size $k$
appears in exactly one of $b_k$ possible types. For more details
on properties of multisets, we refer the reader e.g. to [3].
Weighted partitions are also associated with the generalized
Bose-Einstein model of ideal gas, where $n(=E)$ is interpreted as
the total energy of the system of particles. The weights $b_k,
k\ge 1$, are viewed as counts of the distinct positions of the
particles in the state space, where a particle in a given position
has (rescaled) energy $k$ (for more details on the relationship
between combinatorial partitions and various models of ideal gas,
see [21]). From combinatorial point of view, it is fairly natural
to assume that $b_k, k\ge 1$, are integers (see e.g. the``money
changing problem'' discussed in detail in [24; Sect. 3.15]). On
the other hand, it turns out that this requirement is not
necessary for the analytical approach used in this paper. That is
why, we assume that $b_k, k\ge 1$, are real non-negative numbers.

For a given sequence $b=\{b_k, k\ge 1\}$, let $\mathcal{P}_b(n)$
be the set of all weighted partitions of the positive integer $n$
and let $p_b(n)=\mid\mathcal{P}_b(n)\mid$ be its cardinality. It
is known that the generating function $f_b(x)$ of the numbers
$p_b(n)$ is of Euler's type, namely,
\begin{equation}\label{euler}
f_b(x)=1+\sum_{n=1}^\infty p_b(n)x^n =\prod_{k=1}^\infty
(1-x^k)^{-b_k},\quad \mid x\mid<1
\end{equation}
(see [24; Sect. 3.14]).
 We introduce the uniform probability measure
 $\mathbb{P}=\mathbb{P}_{n,b}$ on the set of weighted partitions of $n$ assuming that the
 probability $1/p_b(n)$ is assigned to each $n$-partition with weight sequence $b$.
 In this paper we focus on the size of the largest part $X_n$ of a
 random weighted partition of $n$. With respect to the probability
 measure $\mathbb{P}$, $X_n$ becomes a random variable, defined on the set
 $\mathcal{P}_b(n)$. It is also well-known that
\begin{equation}\label{df}
f_{m,b}(x)=1+\sum_{n=1}^\infty p_b(n)\mathbb{P}(X_n\le m)x^n
=\prod_{k=1}^m (1-x^k)^{-b_k},\quad m\ge 1
\end{equation}
(see [24; Sect. 3.15]).

 The asymptotic behavior of the combinatorial numbers $p_b(n)$
 (the Taylor coefficients in (\ref{euler})) will play important role in our further analysis. A
 fairly general scheme of assumptions on the parametric
 sequence $b$ was proposed by Meinardus [13] (see also [2; Chap. 6]), who
 found
 an asymptotic expansion for the numbers $p_b(n)$ as $n\to\infty$.
 The same asymptotic was also obtained in [10], where one of the
 Meinardus conditions was weakened.
 Meinardus' approach
 is based on considering two generating series:
 \begin{equation}
 D(s)=\sum_{k=1}^\infty b_k k^{-s}, \quad s=\sigma+iy, \label{dofes}
 \end{equation}
 and
 \begin{equation}
 G(z)=\sum_{k=1}^\infty b_k z^k,\quad \mid z\mid\le 1. \label{gfbk}
 \end{equation}
 Below we give the Meinardus scheme of conditions. Throughout the
 paper by $\Re(z)$ and $\Im(z)$ we denote the real and imaginary part
 of the complex number
 $z$, respectively.

$(M_1)$ The Dirichlet series (\ref{dofes}) converges in the
half-plane $\sigma>\rho>0$ and there is a constant $C_0\in(0,1)$,
such that the function $D(s)$ has an analytic continuation to the
half-plane $\{s:\sigma\ge -C_0\}$ on which it is analytic except
for the simple pole at $s=\rho$ with residue $A>0$.

$(M_2)$ There exists a constant $C_1>0$ such that
 $$
 D(s)=O(\mid y\mid^{C_1}),\quad \mid y\mid\to\infty,
 $$
 uniformly for $\sigma\ge -C_0$.

 $(M_3)$ There are constants $\epsilon>0$ and $C_2(=C_2(\epsilon)>0)$, such that the
 function $g(\tau)=G(e^{-\tau}), \tau=\alpha+2\pi iu$, $u$ real
 and $\alpha>0$ (see (\ref{gfbk})) satisfies
 $$
 \Re(g(\tau))-g(\alpha)\le-C_2\alpha^{-\epsilon},\quad \mid\arg(\tau)\mid>\pi/4,
 \quad 0\neq\mid u\mid\le 1/2,
 $$
 for enough small values of $\alpha$.

The first assumption (M1) specifies the
 domain, say $\mathcal{H}$, in which $D(s)$ has an analytic
 continuation. The second is related to the asymptotic behavior of
 $D(s)$, whenever $|\Im(s)|\to\infty$. Functions, which are bounded
 by $O(|\Im(s)|^q), 0<q<\infty$, in certain domain, as
 $|\Im(s)|\to\infty$,
 are called
 functions of finite order. It is known that the sum of the Dirichlet
 series in (\ref{dofes}) satisfies the finite order property in
 a closed half-plane contained in the
 half-plane of convergence $\sigma>\rho$ (see e.g. [20; Sect. 9.3.3]).
 The Meinardus second condition requires that the same holds for the
 analytic continuation of $D(s)$ in the whole domain $\mathcal{H}$. Finally,
 the Meinardus third condition implies a bound on
 $\Re(G(e^{-\tau}))$ (see (\ref{gfbk})) for certain specific complex
 values of $\tau$. In some cases its verification is technically
complicated. Since
\begin{equation}\label{regtau}
\Re(g(\tau))-g(\alpha) =-2\sum_{k=1}^\infty b_k
e^{-k\alpha}\sin^2{(\pi ku)}
\end{equation}
and the inequality $\mid\arg(\tau)\mid>\pi/4$ implies that
$\tan{(\mid\arg(\tau)\mid)}=2\pi\mid u\mid/\alpha>1$, condition
$(M_3)$ can be also reformulated as follows:
 \begin{equation}\label{sinus}
 S_n:=\sum_{k=1}^\infty b_k e^{-k\alpha}\sin^2{(\pi ku)}
 \ge C_2\alpha^{-\epsilon}, \quad 0<\frac{\alpha}{2\pi}<\mid
 u\mid\le 1/2,
 \end{equation}
 for small enough $\alpha$ and some constants $C_2, \epsilon>0$ $(C_2=C_2(\epsilon))$
 [10; p.310].

 Moreover, Granovsky et al. [10; Lemma 1] proved that this inequality holds for any sequence
 $b_k, k\ge 1$, satisfying the inequality $b_k\ge Ck^{\nu-1}, k\ge
 k_0$, for some $k_0\ge 1$ and $C, \nu>0$. We notice that if
 \begin{equation}\label{cknu}
 b_k=C k^{\nu-1}, \quad k\ge 1,
 \end{equation}
 then $D(s)=C\zeta(s-\nu+1)$, where $\zeta$ denotes the Riemann
 zeta function. Therefore, $D(s)$ has a single pole at $s=\nu$
 with residue $C>0$ and a meromorphic analytic continuation to the
 whole complex plane [23; Sect. 13.13]. These facts show that
 conditions $(M_1)-(M_3)$ are satisfied by the weights
 (\ref{cknu}) with $\rho=\nu$ and $A=C$.

 Throughout the paper we assume that conditions $(M_1)-(M_3)$
  are satisfied. Our aim is to determine asymptotically, as
 $n\to\infty$, the distribution of the maximal part size $X_n$.
 Recalling (\ref{df}), we also point out that our results may be interpreted
 in terms of the asymptotic of the combinatorial counts
 of partitions whose part sizes are $\le m$, where the range of the values of $m$ is specified by the weak convergence
 of the random variable $X_n$ to a non-degenerate
 random variable. In the brief review
 given below we summarize some known results on the limiting
 behavior of the random variable $X_n$.

 Consider first the classical case of linear integer partitions,
 where the weights satisfy $b_k=1, k\ge 1$. This kind of
 partitions were broadly studied by many authors in many respects.
 Their graphical representations by Ferrers diagrams show that
 their total number of parts and their maximal part size $X_n$ are
 identically distributed for all $n$ (see [2; Sect. 1.3]). Erd\"{o}s and Lehner [6] were
 apparently the first who applied a probabilistic approach to the
 study of integer partitions.
 As a matter of fact, they found
 an appropriate normalization for $X_n$ in this case and showed that
 $\pi X_n/(6n)^{1/2}-\log{((6n)^{1/2}/\pi)}$ converges weakly, as $n\to\infty$, to
 a random variable having the extreme value (Gumbel)
 distribution. The local version of their theorem was derived later by Auluck et al. [4].
 Fristedt [9] studied linear integer partitions using a transfer
 method to functionals of independent and geometrically
 distributed random variables. Among other results, he obtained
 the limiting distribution of the $k$th largest part size whenever
 $k$ is fixed. Finally, we notice that among weighted integer partitions
 only the linear ones possess the property that the number of parts and the maximum part size
 are identically distributed.
 The limiting distribution of the number of parts in the general case of
 random weighted
 partitions under the Meinardus scheme of conditions is studied in
 [15]. It turns out that the limiting distribution laws depend
 on particular ranges in which the parameter $\rho$ varies (see
 condition $(M_1)$).

Another important particular case of weighted partitions arises
whenever $b_k=k, k\ge 1$. It turns out that in this case the
generating function $f_b(x)$ (see (\ref{euler})) enumerates the
plane partitions. A plane partition of $n\ge 1$ is a matrix of
non-negative integers arranged in non-increasing order from left
to right and from top to bottom, so that the double sum of its
elements equals $n$. Together with the largest part size $X_n$,
consider also the counts of the non-zero rows and columns of the
matrix of a plane partition. It turns out that these three
quantities measure the sizes of the corresponding solid diagram of
a plane partition in the $3D$ space. (The solid diagram is a heap
of $n$ unit cubes placed in the first octant of a coordinate
system in a $3D$ space whose columns composed by stacked cubes
have non-increasing heights along the $x$- and $y$-axis; the
height of this heap along the $z$-axis is just $X_n$, the largest
part size.) Similarly to Ferrers diagrams for linear integer
partitions, the three sizes of this heap appear to be identically
distributed for all $n\ge 1$ (for more details, see [19; p. 371]).
Their joint limiting distribution was found in [17]. The marginal
limiting distributions (including the limiting distribution of
$X_n$) were obtained in [14]. For more details on various
properties of plane partitions and their applications to
combinatorics and analysis of algorithms, we refer the reader to
[2; Chap. 11], [16; Chap. 11] and [19; Chap. 7].

Our study is also closely related to some recent results on the
maximal particle energy in the Bose-Einstein model of ideal gas.
The general setting and the probabilistic frame of problems from
statistical mechanics and their relationship with enumerative
combinatorics were given by Vershik [21]. In the context of the
infinite product formula (\ref{euler}), he studied the
Bose-Einstein model by a family of probability measures $\mu_v,
v\in (0,1)$, defined on the set of all $b$-weighted partitions
$\mathcal{P}_b =\cup_{n\ge 0}\mathcal{P}_b(n)$. So, for a
partition $\lambda=(\lambda_1,...,\lambda_l)\in\mathcal{P},
\lambda_1\ge...\ge\lambda_l>0$, let
$r_k(\lambda)=\{j:\lambda_j=k\}$ denote the number of parts of
$\lambda$ that are equal to $k\ge 1$. Then $\mu_v$ is defined by
$$
\mu_v(\{\lambda\in\mathcal{P}: r_k(\lambda)=j\})={b_k+j-1\choose
j}
 v^{kj}(1-v^k)^{b_k}, \quad 0<v<1.
 $$
 The key feature in the study of this kind of measures is the fact
 that a kind of a conditional probability measure on
 $\mathcal{P}_b(n)$ turns out to be independent of $v$ for all $n$
 and coincides with the uniform probability measure $\mathbb{P}=\mathbb{P}_{n,b}$
 (for more details, see [21]). In [22] Vershik
 and Yakubovich studied
  the limiting distribution of the maximal
 particle energy $X(\lambda)$, or, which is the same, the largest part
 size $X(\lambda), \lambda\in\mathcal{P}$, with respect to the measure $\mu_v$, as $v\to
 1^-$. In particular, under the assumption that the weight sequence $b$
 satisfies (\ref{cknu}), they proved that
 \begin {eqnarray}\label{veya}
 & & \lim_{v\to 1^-}\mu_v(\{\lambda\in\mathcal{P}: (1-v)X(\lambda)
 -\nu\mid\log{(1-v)}\mid \nonumber \\
 & & -(\nu-1)\log{\mid\log{(1-v)}\mid}
 -(\nu-1)\log{\nu}-\log{C}\le t\}) \nonumber \\
 & & =e^{-e^{-t}}, \quad -\infty<t<\infty.
 \end{eqnarray}
 As it was mentioned before, the weight sequence (\ref{cknu})
 satisfies conditions $(M_1)-(M_3)$. Vershik and
 Yakubovich [22] studied also a more realistic model of quantum ideal
 gas in a $N$-dimensional space, for which the weights satisfy $\sum_{j=1}^k
 b_j=K_N k^{N/2}+O(k^{\kappa_N})$, as $k\to\infty$ ($K_N$ and $\kappa_N<N/2$ are computable constants).

The main result of this paper is obtained in terms of the uniform
probability measure $\mathbb{P}=\mathbb{P}_{n,b}$ on the set
$\mathcal{P}_b(n)$. Before stating it, for the sake of brevity, we
introduce the following notation:
\begin{equation}
 a(n)=a(n;\rho,A)
 =\left(\frac{A\Gamma(\rho+1)\zeta(\rho+1)}{n}\right)^{\frac{1}{\rho+1}},
 \quad n\ge 1,
 \label{norm}
 \end{equation}
 where the constants $\rho$ and $A$ are defined by condition
 $(M1)$.

 \begin{theorem} If the weight sequence $b$, satisfies
 conditions $(M1)-(M3)$, then, for all real $t$, the limiting
 distribution of the largest part size $X_n$ is given by
\begin{eqnarray}
 & & \lim_{n\to\infty} \mathbb{P}(a(n)X_n+\rho\log{a(n)}
 -(\rho-1)\log{\mid\log{a(n)}\mid}-(\rho-1)\log{\rho}-\log{A} \le
 t) \nonumber \\
 & & =e^{-e^{-t}}.\label{limth}
 \end{eqnarray}
 \end{theorem}

{\it Remark.} One can easily compare (\ref{limth}) with
(\ref{veya}) setting in the latter one $v=1-a(n), \nu=\rho$ and
$C=A$ and observe the coinciding normalizations. We also notice
that the limiting results for linear and plane partitions (see
[6,14]) follow from (\ref{limth}) with $A=1$ and $\rho=1$ and $2$,
respectively.

 The method of our proof combines Hayman's theorem for estimating
 coefficients of admissible power series [11] (see also [7; Sect.
 VIII.5]), a generalization of Perron formula, which yields the expression for partial
 sums of a Dirichlet series by a complex integral of the inverse
 Mellin transform applied to the Dirichlet series itself (see
 Thm 3.1 from the Supplement of [18]) and some Mellin transform computations.

 We organize our paper as follows. Section 2 includes
  some auxiliary
  facts that we need further. Some proofs are
 omitted since they are given in [10,13]. In Section 3 we present the proof of Theorem 1.
 The Appendix contains some technical details related to the
 application of the generalized Perron formula [18].

\section{Preliminary Results}
\setcounter{equation}{0}

 We start with a lemma establishing an asymptotic estimate for
 infinite product representation (\ref{euler}) of the
 generating function $f_b(x)$. It has been proved by Meinardus
 [13] (see also [2; Lemma 6.1]).

 \begin{lemma}
  Suppose that sequence $b$ is such that the
 associated Dirichlet series (\ref{dofes}) satisfies conditions
 $(M_1)$ and $(M_2)$. If $\tau=\alpha+i\theta$, then
 $$
 f_b(e^{-\tau})
 =\exp{(A\Gamma(\rho)\zeta(\rho+1)\tau^{-\rho}-D(0)\log{\tau}
 +D^\prime(0)+O(\alpha^{C_0}))}
 $$
 as $\alpha\to 0^+$ uniformly for $\mid\theta\mid\le\pi$ and
 $\mid\arg{\tau}\mid\le\pi/4$. $\rule{2mm}{2mm}$
 \end{lemma}

 Our first goal is to show that Meinardus' conditions $(M_1)-(M_3)$
 imply that the generating function $f_b(x)$ possesses the Hayman
 admissibility properties [11] (see also [7; Sect. VIII.5]) in the
 unit disc. For $0<r<1$, we introduce the
 functions:
 \begin{equation}\label{capitalf}
 F_b(r) =\log{f_b(r)} =-\sum_{k=1}^\infty b_k\log{(1-r^k)},
 \end{equation}
 \begin{equation}\label{ar}
 \mathcal{A}_b(r) =rF_b^\prime(r) =r\frac{f_b^\prime(r)}{f_b(r)},
 \end{equation}
 \begin{equation}\label{br}
 \mathcal{B}_b(r) =r^2 F_b^{\prime\prime}(r)+rF_b^{\prime}(r) =r\frac{f_b^\prime(r)}{f_b(r)}
 +r^2\frac{f_b^{\prime\prime}(r)}{f_b(r)}
 -r^2\left(\frac{f_b^\prime(r)}{f_b(r)}\right)^2.
 \end{equation}
 Furthermore, setting in (\ref{capitalf})-(\ref{br})
 $r=e^{-\alpha}$, we shall obtain their asymptotic expansions as
 $\alpha\to 0^+$. For the sake of convenience, we also set
 \begin{equation}\label{consth}
 h=h(\rho,A)=A\Gamma(\rho+1)\zeta(\rho+1).
 \end{equation}
 The proof of the next lemma is contained in [10; Lemma 2].
 \begin{lemma}
  Conditions $(M_1)$ and $(M_2)$ imply the
 following asymptotic expansions:
 \begin{equation}\label{aalpha}
 \mathcal{A}_b(e^{-\alpha}) =h\alpha^{-\rho-1} +D(0)\alpha^{-1}
 +D^\prime(0)+O(\alpha^{C_0-1}),
 \end{equation}
 \begin{equation}\label{balpha}
 \mathcal{B}_b(e^{-\alpha})
 =\frac{d}{d\alpha}(-\mathcal{A}_b(e^{-\alpha}))
 =h(\rho+1)\alpha^{-\rho-2} +D(0)\alpha^{-2} +O(\alpha^{C_0-2}),
 \end{equation}
 \begin{equation}\label{terz}
 F_b^{\prime\prime\prime}(e^{-\alpha}) =O(\alpha^{-\rho-3}),
 \end{equation}
 as
 $\alpha\to 0^+$, where $F_b , \mathcal{A}_b$ and $\mathcal{B}_b$ are
  defined by (\ref{capitalf})-(\ref{br}), respectively. Moreover,
 from (\ref{aalpha}) it follows that the equation
 \begin{equation}\label{eqhay}
 \mathcal{A}_b(e^{-\alpha})=n, \quad n\ge 1,
 \end{equation}
 has a unique solution $\alpha=\alpha_n$, such that $\alpha_n\to 0$
 as $n\to\infty$. An asymptotic expansion of this solution, as
 $n\to\infty$, is given by
 \begin{equation}\label{alphan}
 \alpha_n=a(n) +\frac{D(0)}{(\rho+1)n} +O(n^{-1-\beta}),
 \end{equation}
 where
 $\beta=\min{\left(\frac{C_0}{\rho+1},\frac{\rho}{\rho+1}\right)}$
 and $a(n)$ are the normalizing constants given by (\ref{norm}). $\rule{2mm}{2mm}$
 \end{lemma}

 We notice that (\ref{aalpha}), (\ref{balpha}) and (\ref{alphan})
 imply that
 \begin{equation} \label{capture}
 \mathcal{A}_b(e^{-\alpha_n})\to\infty, \quad
 \mathcal{B}_b(e^{-\alpha_n})\to\infty, \quad n\to\infty,
 \end{equation}
 that is, Hayman ``capture'' condition [7; p. 565] is satisfied
 with
 $r=r_n=e^{-\alpha_n}$. Our next step is to establish Hayman
  ``locality'' condition, which implies the asymptotic behavior of
 $f_b(x)$ in a suitable neighborhood of $x=1$.

 \begin{lemma} Suppose that the weight sequence $b$ satisfies conditions $(M_1)$ and
 $(M_2)$ and $\alpha_n$ is the solution of (\ref{eqhay}) given by
 (\ref{alphan}).  Let
 \begin{equation} \label{delta}
 \delta_n=\alpha_n^{1+\rho/3}/\omega(n), \quad n\ge 1,
 \end{equation}
 where $\omega(n)\to\infty$ as $n\to\infty$ arbitrarily slowly.
 Then
 \begin{equation}\label{locality}
 e^{-i\theta
 n}\frac{f_b(e^{-\alpha_n+i\theta})}{f_b(e^{-\alpha_n})}
 =e^{-\theta^2 \mathcal{B}_b(e^{-\alpha_n})/2}(1+O(1/\omega^3(n))
 \end{equation}
 uniformly for $\mid\theta\mid\le\delta_n$.
 \end{lemma}

{\it Proof.} Applying Lemma 1, we observe that
\begin{eqnarray}\label{malko}
& & e^{-i\theta
 n}\frac{f_b(e^{-\alpha_n+i\theta})}{f_b(e^{-\alpha_n})}
 \\
 & & =\exp{\left(\frac{h}{\rho}((\alpha_n-i\theta)^{-\rho}
 -\alpha_n^{-\rho})-D(0)\log{\left(1-\frac{i\theta}{\alpha_n}\right)}
 -i\theta n +O(\alpha_n^{C_0})\right)}, \nonumber
\ \end{eqnarray}
 where $h$ is given by (\ref{consth}). Expanding
 $(\alpha_n-i\theta)^{-\rho}$ and $\log{(1-i\theta/\alpha_n)}$ by
 Taylor formula and using (\ref{aalpha}), (\ref{balpha}) and
 (\ref{eqhay}), we obtain
 \begin{eqnarray}
 & & \frac{h}{\rho}((\alpha_n-i\theta)^{-\rho}-\alpha_n^{-\rho})
 -D(0)\log{\left(1-\frac{i\theta}{\alpha_n}\right)} -i\theta n
 \nonumber \\
  & & =i\theta(h\alpha_n^{-\rho-1} +D(0)\alpha_n^{-1} +D^\prime(0)-n)
 -\frac{\theta^2}{2} h(\rho+1)\alpha_n^{-\rho-2}
 -\frac{D(0)\theta^2}{2\alpha_n^2} \nonumber \\
 & & -i\theta D^\prime(0)
  +O(\mid\theta\mid^3\alpha_n^{-3-\rho})
 =i\theta(\mathcal{A}_b(e^{-\alpha_n}) -n +O(\alpha_n^{C_0-1}))
 \nonumber \\
 & & -\frac{\theta^2}{2}(\mathcal{B}_b(e^{-\alpha_n})
  +O(\alpha_n^{C_0-2}))
  +O(\mid\theta\mid) +O(\mid\theta\mid^3\alpha_n^{-3-\rho})
 \nonumber \\
 & & =-\frac{\theta^2}{2}\mathcal{B}_b(e^{-\alpha_n})
  +O(\delta_n\alpha_n^{C_0-1}) +O(\delta_n^2\alpha_n^{C_0-2})
 +O(\delta_n)+O(\delta_n^3\alpha_n^{-3-\rho}). \nonumber
 \end{eqnarray}
 Substituting this into  (\ref{malko}) and taking into account
 (\ref{delta}), we obtain (\ref{locality}). $\rule{2mm}{2mm}$

  To study the
 behavior of $f_b(e^{-\alpha_n+i\theta})$ outside the range
 $-\delta_n<\theta<\delta_n$ we need the Wiener-Ikehara Tauberian theorem
 on Dirichlet series. It tells us how condition $(M_1)$ implies
 an asymptotic estimate for the partial sums $\sum_{k=1}^n b_k$.

 {\bf Wiener-Ikehara Theorem.} {\it (See [12; Thm. 2.2, p. 122].)
 Suppose that the Dirichlet series $\tilde{D}(s)=\sum_{k=1}^\infty
 c_k k^{-s}$ is such that the function
 $\tilde{D}(s)-\frac{c}{s-1}$ has an analytic continuation to the
 closed half-plane $\Re(s)\ge 1$. Then
 \begin{equation}\label{wi}
 \sum_{k=1}^n c_k\sim cn, \quad n\to\infty. \quad \rule{2mm}{2mm}
 \end{equation}}

 We also denote by $\{\gamma\}$ the fractional part of the
 real number $\gamma$, and by $\parallel\gamma\parallel$ the
 distance from $\gamma$ to the nearest integer, so that
 \begin{equation} \label{gam}
 \parallel\gamma\parallel=\left\{\begin{array}{ll} \{\gamma\} & \qquad  \mbox {if}\qquad \{\gamma\}\le 1/2, \\
 1-\{\gamma\} & \qquad \mbox {if}\qquad \{\gamma\}>1/2.
 \end{array}\right.
 \end{equation}
 It is not difficult to show that
 \begin{equation} \label{sinlowbound}
 \sin^2{(\pi\gamma)}\ge 4\parallel\gamma\parallel^2
 \end{equation}
 (see [8; p. 272]). Now, we are ready to prove that Hayman
 last (``decay'') condition [7; p. 565] is also valid.

 \begin{lemma} Suppose that $f_b(x)$ satisfies
 conditions $(M_1)-(M_3)$. Then, for sufficiently large $n$,
 $$
 \mid f_b(e^{-\alpha_n+i\theta})\mid\le f_b(e^{-\alpha_n})
 e^{-C_3\alpha_n^{-\epsilon_1}}
 $$
 uniformly for $\delta_n\le\mid\theta\mid<\pi$, where $C_3$ and $\epsilon_1$ are
 positive
 constants.
 \end{lemma}

{\it Proof.} First, we notice that
\begin{equation} \label{modul}
\frac{\mid f_b(e^{-\alpha_n+i\theta})\mid}{f_b(e^{-\alpha_n})}
=\exp{(\Re({\log{f_b(e^{-\alpha_n+i\theta}))}}
-\log{f_b(e^{-\alpha_n})})}.
\end{equation}
  Then, setting $\theta=2\pi u$, for
 $\alpha_n/2\pi\le\mid u\mid=\mid\theta\mid/2\pi<1/2$, we almost repeat the argument
  from [10; p. 324]:
  \begin{eqnarray} \label{realparts}
  & & \Re{(\log{f_b(e^{-\alpha_n+i\theta})})}
-\log{f_b(e^{-\alpha_n})} \nonumber \\
 & & =\Re{\left(-\sum_{k=1}^\infty b_k \log{\left(\frac{1-e^{-k\alpha_n+2\pi
 iuk}}{1-e^{-k\alpha_n}}\right)}\right)} \nonumber \\
  & & =-\frac{1}{2}\sum_{k=1}^\infty b_k
  \log{\left(\frac{1-2e^{-k\alpha_n}\cos{(2\pi uk)}+e^{-2\alpha_n k}}
  {(1-e^{-\alpha_n k})^2}\right)} \nonumber \\
  & & =-\frac{1}{2}\sum_{k=1}^\infty b_k
  \log{\left(1+\frac{4e^{-\alpha_n k}\sin^2{(\pi uk)}}{(1-e^{-\alpha_n
  k})^2}\right)} \nonumber \\
  & & \le -\frac{1}{2}\sum_{k=1}^\infty b_k \log{(1+4e^{-\alpha_n k}\sin^2{(\pi
  uk)})} \nonumber \\
  & & \le -\frac{\log{5}}{2}\sum_{k=1}^\infty b_k e^{-\alpha_n
  k}\sin^2{(\pi uk)} \nonumber \\
  & & =-\frac{\log{5}}{2} S_n \le -\frac{\log{5}}{2}
  C_2\alpha_n^{-\epsilon},
 \end{eqnarray}
 where the last two inequalities follow from the fact that
 $\log{(1+y)}\ge\left(\frac{\log{5}}{4}\right)y$ $(0\le y\le 4)$ and (\ref{sinus}), respectively.
 Thus, the required inequality is proved for $\alpha_n\le\mid\theta\mid<\pi$.
 It remains to consider the interval $\delta_n\le\mid\theta\mid<\alpha_n$.
 (\ref{regtau}) implies that
 we have to find now a lower bound for the sum $S_n$ in (\ref{sinus}) if $\delta_n/2\pi\le\mid u\mid<\alpha_n/2\pi$
 (i.e. if $\delta_n\le\mid\theta\mid<\alpha_n$).
 We shall apply
 Wiener-Ikehara Tauberian theorem setting there
 $c_k=k^{-\rho+1}b_k, k\ge 1$, and
 \begin{equation}\label{dtil}
 \tilde{D}(s)=\sum_{k\ge 1}
 b_k k^{-s-\rho+1}=D(s+\rho-1), \quad s=\sigma+iy.
 \end{equation}
 Since both $C_0,\rho>0$ from
 condition $(M_1)$, the function $\tilde{D}(s)$ satisfies the
 condition of Wiener-Ikehara theorem with $c=A$. Moreover, since by (\ref{alphan})
 $\alpha_n^{-1}\to\infty$ as $n\to\infty$, we can
 apply (\ref{wi}) in the form
 \begin{equation} \label{wik}
 \sum_{1\le k\le K^/\alpha_n} k^{-\rho+1}b_k\sim AK/\alpha_n,
 \quad n\to\infty,
 \end{equation}
 where the constant $K>0$ will be specified later. Our next argument will be similar
 to that given in [8; Lemma 7]. First, using (\ref{gam}), we
 observe that $\parallel uk\parallel=uk$ if $\mid
 u\mid k<1/2$. This implies that, for $1\le k\le\pi/\alpha_n$, $\parallel
 uk\parallel$ can be replaced by $\mid u\mid k$. Recalling that
 $\mid u\mid\ge\delta_n/2\pi$ and applying (\ref{sinlowbound}) and
 (\ref{delta}), we obtain
 \begin{eqnarray}
 & & S_n =\sum_{k=1}^\infty b_k e^{-k\alpha_n}\sin^2{(\pi uk)}
 \ge 4\sum_{k=1}^\infty b_k e^{-k\alpha_n}\parallel uk\parallel^2
 \nonumber \\
 & & \ge 4u^2\sum_{1\le k\le\pi/\alpha_n} k^2 b_k e^{-k\alpha_n}
 \ge (\delta_n/\pi)^2 \sum_{1\le k\le\pi/\alpha_n} k^2 b_k
 e^{-k\alpha_n} \nonumber \\
 & & =\frac{\alpha_n^{2+2\rho/3}}{(\pi\omega(n))^2}
 \sum_{1\le k\le\pi/\alpha_n} k^{\rho+1}\frac{k^2 b_k}{k^{\rho+1}}
 e^{-k\alpha_n} \nonumber \\
 & & =\frac{\alpha_n^{2+2\rho/3}}{(\pi\omega(n))^2} \sum_{1\le
 k\le\pi/\alpha_n} (k^{\rho+1} e^{-k\alpha_n})(b_k k^{-\rho+1}).
 \nonumber
 \end{eqnarray}
 It is easy to check that the sequence $k^{\rho+1}e^{-k\alpha_n},
 k\ge 1$, is non-increasing if
 $k\ge(\rho+1)/\alpha_n-1/2+O(\alpha_n)$. Hence, for $\rho+1<\pi$,
 (\ref{wik}) with $K=\pi,\rho+1$ implies that
 \begin{eqnarray} \label{snpi}
 & & S_n \ge\frac{\alpha_n^{2+2\rho/3}}{(\pi\omega(n))^2} \sum_{(\rho+1)/\alpha_n\le
 k\le\pi/\alpha_n} (k^{\rho+1} e^{-k\alpha_n})(b_k k^{-\rho+1})
 \nonumber \\
 & & \ge\frac{\alpha_n^{2+2\rho/3}}{(\pi\omega(n))^2}
 (\pi/\alpha_n)^{\rho+1} e^{-\pi}\sum_{(\rho+1)/\alpha_n\le
 k\le\pi/\alpha_n} b_k k^{-\rho+1} \nonumber \\
 & & =\frac{\pi^{\rho-1}e^{-\pi}\alpha_n^{1-\rho/3}}
 {\omega^2(n)}\left(\sum_{1\le k\le\pi/\alpha_n} b_k k^{-\rho+1}
 -\sum_{1\le k<(\rho+1)/\alpha_n} b_k k^{-\rho+1}\right) \nonumber \\
 & & \sim A\pi^{\rho-1}(\pi-\rho-1)e^{-\pi}\alpha_n^{-\rho/3}/
 \omega^2(n).
 \end{eqnarray}
 If $\pi\le\rho+1$, then $k^{\rho+1}e^{-k\alpha_n}$ is a
 non-decreasing sequence for $1\le k\le\pi/\alpha_n$. Then, for
 some $l\in (0,\pi)$, $k^{\rho+1} e^{-k\alpha_n}\ge l^{\rho+1}
 e^{-l\alpha_n}$ and in the same way we observe that
 \begin{eqnarray}\label{snro}
 & & S_n\ge \frac{\alpha_n^{2+2\rho/3}}{(\pi\omega(n))^2} \sum_{l/\alpha_n\le
 k\le\pi/\alpha_n} (k^{\rho+1} e^{-k\alpha_n})(b_k k^{-\rho+1})
 \nonumber \\
 & & \ge \frac{\alpha_n^{2+2\rho/3} l^{\rho+1} e^{-l}}{(\pi\omega(n))^2 \alpha_n^{\rho+1}}
 \sum_{l/\alpha_n\le k\le\pi/\alpha_n} b_k k^{-\rho+1}
  \nonumber \\
  & & \sim
 \frac{Al^{\rho+1}e^{-l}(\pi-l)\alpha_n^{-\rho/3}}
 {(\pi\omega(n))^2}(1+o(1)).
 \end{eqnarray}
 Consequently, (\ref{realparts}), (\ref{snpi}) and
 (\ref{snro}) imply that there are two constants $C_3,\epsilon_1>0$, such that
 \begin{equation}\label{sn}
 S_n \ge C_3\alpha_n^{-\epsilon_1}
 \end{equation}
 uniformly for $\delta_n/2\pi\le\mid u\mid<1/2$.
 Moreover, $\epsilon_1\le\min{(\epsilon,\rho/3)}$ since $\omega(n)\to\infty$ as $n\to\infty$
 arbitrarily slowly.
 Hence, noting that $u=\theta/2\pi$, we obtain that the required inequality
 holds uniformly for $\delta_n\le\mid\theta\mid<\pi$.  \rule{2mm}{2mm}

 We now recall (\ref{balpha}) from Lemma 2. It implies that
 $$
 \mathcal{B}_b^{1/2}(e^{-\alpha_n})\sim
 (h(\rho+1))^{1/2}\alpha_n^{-1-\rho/2}, \quad n\to\infty.
 $$
 Combining this asymptotic equivalence with the result of Lemma 4,
 we obtain Hayman ``decay'' condition [7; p. 565], namely,
\begin{equation}\label{decay}
 \mid f_b(e^{-\alpha_n+i\theta})\mid
 =o(f_b(e^{-\alpha_n})/\mathcal{B}_b^{1/2}(e^{-\alpha_n})), \quad
 n\to\infty,
 \end{equation}
 uniformly for $\delta_n\le\mid\theta\mid<\pi$.

 Eqs. (\ref{capture}), (\ref{locality}) and (\ref{decay}) show
 that the function $f_b(x)$ is admissible in the sense of Hayman.
 Therefore, we can apply Thm. VIII.4 of [7] for its coefficients.
 We state this result in the next lemma.

 \begin{lemma} Suppose that the weight sequence $b$ satisfies
 Meinardus conditions $(M_1)-(M_3)$. Then, the asymptotic for the
 total number of weighted partitions is given by
 \begin{equation} \label{totalnr}
 p_b(n) \sim \frac{e^{n\alpha_n} f_b(e^{-\alpha_n})}{\sqrt{2\pi\mathcal{B}_b(e^{-\alpha_n})}}
 \end{equation}
 as $n\to\infty$, where $\alpha_n$ is the unique solution of
 (\ref{eqhay}) whose asymptotic expansion is given by
 (\ref{alphan}) and $\mathcal{B}_b(e^{-\alpha_n})$ is defined by
 (\ref{balpha}). $\rule{2mm}{2mm}$
 \end{lemma}

 {\it Remark.} The asymptotic equivalence (\ref{totalnr}) is in
 fact Meinardus asymptotic formula [13] for the number of weighted
 partitions of $n$. Here we give the formula in
 a slightly different form, which is more convenient for our further
 asymptotic analysis. One can easily show the coincidence of
 (\ref{totalnr}) with the Meinardus original formula, applying the
 result of Lemma 1 to $f_b(e^{-\alpha_n})$ and replacing
 $\alpha_n$ and $\mathcal{B}_b(e^{-\alpha_n})$ by
 (\ref{alphan}) and (\ref{balpha}), respectively.

 Further, we also need a bound on the rate of growth of the
 weights $b_k$, as $k\to\infty$.
 Using Wiener-Ikehara
 Tauberian theorem, Granovsky et al. [10; p. 310] showed that
 $b_k=o(k^\rho)$ as $k\to\infty$. We need this bound in a
 different slightly more precise form.

 \begin{lemma} If the sequence of weights $b$
 satisfies conditions $(M_1)$ and $(M_2)$, then there is a
 sequence of numbers $L_k, k\ge 1$, satisfying
 $\lim_{k\to\infty}L_k=0$ and such that
 \begin{equation}\label{recur}
 b_k=(L_k-L_{k-1})k^\rho+(A+L_{k-1})k^{\rho-1}, \quad k\ge 2,
 \quad L_1=b_1-A,
 \end{equation}
 where $A$ is the constant defined in condition $(M_1)$.
 \end{lemma}

 {\it Proof.} As in [10; p. 310], we rewrite (\ref{wi}) in the
 following way:
 $$
 \frac{1}{k}\sum_{j=1}^k c_j
 =\frac{1}{k}c_k+\frac{1}{k}\sum_{j=1}^{k-1}c_j =c+L_k, \quad k\ge 2,
 $$
 where $\lim_{k\to\infty}L_k=0$. We also set $L_1=c_1-c$. Then, for $k\ge 2$, we have
 $$
 \frac{1}{k}c_k=L_k-L_{k-1}+\frac{1}{k}(c+L_{k-1}).
 $$
 To obtain
 (\ref{recur}) it remains to set $c_k=k^{-\rho+1}b_k, k\ge 1,
 c=A$ and recall that the weights $b_k$ satisfy conditions
 $(M_1)$ and $(M_2)$ with $C_0, \rho>0$. $\quad\rule{2mm}{2mm}$

Now, we also recall formula (\ref{df}) for the truncated products
$f_{m,b}(x), m\ge 1$. Similarly to (\ref{capitalf}), we set
\begin{equation} \label{capitalfm}
F_{m,b}(x) =\log{f_{m,b}(x)} =-\sum_{k=1}^m b_k\log{(1-x^k)}.
\end{equation}
(Here we consider the main branch of the logarithmic function,
assuming that $\log{y}<0$ for $0<y<1$). Further on, when computing
the derivatives of (\ref{capitalf}) and (\ref{capitalfm}), we
shall write
$$
F_b^{(j)}(e^{-\alpha_n}) =F_b^{(j)}(x)\mid_{x=e^{-\alpha_n}},
\quad F_{m,b}^{(j)}(e^{-\alpha_n})
=F_{m,b}^{(j)}(x)\mid_{x=e^{-\alpha_n}}, j=1,2,3.
$$
Our next lemma establishes estimates on the tails
$F_b^{(j)}(e^{-\alpha_n})-F_{m,b}^{(j)}(e^{-\alpha_n})$ for some
specific values of $m$.

\begin{lemma} Suppose that the weight sequence $b$
satisfies conditions $(M_1)$ and $(M_2)$ and that $\alpha_n, n\ge
1$, is defined by eq. (\ref{alphan}). Moreover, let $m=m(n)$ be a
sequence of integers satisfying
\begin{equation}\label{mn}
m\sim\rho\alpha_n^{-1}\log{\alpha_n^{-1}}, \quad n\to\infty.
\end{equation}
Then
$$
F_b^{(j)}(e^{-\alpha_n})-F_{m,b}^{(j)}(e^{-\alpha_n})
=O(\alpha_n^{-j}\log^{\rho+j}{\alpha_n^{-1}}),\quad j=1,2,3.
$$
\end{lemma}

{\it Proof.} We shall consider only the case $j=1$. The other two cases
are studied in a similar way.

First, we choose a sequence of integers $m_1(n)$ that satisfies
the asymptotic equivalence
\begin{equation}\label{monen}
m_1=m_1(n)\sim (\rho+1)\alpha_n^{-1}\log{\alpha_n^{-1}}
\end{equation}
and decompose the difference of the first derivatives in the
following way:
\begin{equation}\label{diff}
F_b^\prime(e^{-\alpha_n})-F_{m,b}^\prime(e^{-\alpha_n})
=\sum_{k=m+1}^\infty \frac{kb_k
e^{-(k-1)\alpha_n}}{1-e^{-k\alpha_n}} =\Sigma_1+\Sigma_2,
\end{equation}
where
$$
\Sigma_1=\sum_{k=m+1}^{m_1} \frac{kb_k
e^{-(k-1)\alpha_n}}{1-e^{-k\alpha_n}}, \quad \Sigma_2
=\sum_{k=m_1+1}^\infty \frac{kb_k
e^{-(k-1)\alpha_n}}{1-e^{-k\alpha_n}}.
$$
We also notice that (\ref{mn}) and (\ref{monen}) imply that
\begin{equation}\label{expbeh}
e^{-m\alpha_n}\sim\alpha_n^\rho, \quad
e^{-m_1\alpha_n}\sim\alpha_n^{\rho+1}, \quad n\to\infty.
\end{equation}
Hence, applying the result of Lemma 6, for $\Sigma_1$ we obtain
\begin{eqnarray}\label{sigmaprev}
& & \Sigma_1 =O\left(\alpha_n^\rho\sum_{k=m+1}^{m_1}kb_k\right)
\nonumber \\
& &
=O\left(\alpha_n^\rho\sum_{k=m+1}^{m_1}k^{\rho+1}(L_k-L_{k-1})\right)
+O\left(\alpha_n^\rho\sum_{k=m+1}^{m_1}k^\rho\right).
\end{eqnarray}
The second summand in (\ref{sigmaprev}) can be approximated by a
Riemann integral using (\ref{mn}) and (\ref{monen}). We have
$$
\sum_{k=m+1}^{m_1}k^\rho
=m_1^{\rho+1}\sum_{\frac{\rho}{\rho+1}<\frac{k}{m_1}\le 1}
\left(\frac{k}{m_1}\right)^\rho\frac{1}{m_1} \sim m_1^{\rho+1}
\int_{\frac{\rho}{\rho+1}}^1 u^\rho du =O(m_1^{\rho+1}).
$$
Therefore
\begin{equation}\label{secsum}
 O\left(\alpha_n^\rho\sum_{k=m+1}^{m_1}k^\rho\right)
=O(\alpha_n^\rho m_1^{\rho+1})
=O(\alpha_n^{-1}\log^{\rho+1}{\alpha_n^{-1}}).
\end{equation}
The first sum in (\ref{sigmaprev}) can be also rewritten as
\begin{eqnarray}\label{simplif}
& & \sum_{k=m+1}^{m_1}k^{\rho+1}(L_k-L_{k-1}) =-L_m
(m+1)^{\rho+1}+L_{m_1}m_1^{\rho+1} \nonumber \\
& & +\sum_{j=1}^{m_1-m-1}L_{m+j}((m+j)^{\rho+1}
-(m+j+1)^{\rho+1}).
\end{eqnarray}
We recall that by Lemma 6, (\ref{mn}) and (\ref{monen}),
$L_m=o(1)$ and $L_{m_1}=o(1)$. Hence
\begin{equation}\label{t}
(-L_m(m+1)^{\rho+1}+L_{m_1}m_1^{\rho+1})\alpha_n^\rho
=o(\alpha_n^{-1}\log^{\rho+1}{\alpha_n^{-1}}).
\end{equation}
The last sum in (\ref{simplif}) is estimated using again an
approximation by an integral. First, applying a binomial
expansion, we get
\begin{eqnarray}
& & (m+j)^{\rho+1}-(m+j+1)^{\rho+1}
=(m+j)^{\rho+1}\left(1-\left(1+\frac{1}{m+j}\right)^{\rho+1}\right)
\nonumber \\
& & =(m+j)^{\rho+1}\left(-\frac{\rho+1}{m+j}+O((m+j)^{-2})\right)
\nonumber \\
& & =-(\rho+1)(m+j)^\rho+O((m+j)^{\rho-1}). \nonumber
\end{eqnarray}
Then, from (\ref{mn}), (\ref{monen}) and the fact that $L_{m+j}$
are bounded for $1\le j\le m_1-m-1$ we obtain
\begin{eqnarray}\label{bigsum}
& &\sum_{j=1}^{m_1-m-1}L_{m+j}((m+j)^{\rho+1}-(m+j+1)^{\rho+1})
\nonumber\\
& & =-(\rho+1)\sum_{j=1}^{m_1-m-1}
L_{m+j}((m+j)^\rho+O((m+j)^{\rho-1})) \nonumber \\
& & =-(\rho+1)
m_1^{\rho+1}\sum_{\frac{1}{m_1}\le\frac{j}{m_1}<1-\frac{m}{m_1}}
L_{m+j}\left(\frac{1}{m_1}\right)
\left(\left(\frac{m+j}{m_1}\right)^\rho
+O\left(\frac{1}{m_1}\left(\frac{m+j}{m_1}\right)^{\rho-1}\right)\right)
\nonumber \\
& & =O\left(m_1^{\rho+1} \int_0^{1-\frac{\rho}{\rho+1}}
\left(\left(\frac{\rho}{\rho+1}+u\right)^\rho
+O\left(\frac{1}{m_1}\right)\right)du\right) =O(m_1^{\rho+1}).
\end{eqnarray}
Combining (\ref{monen}) with (\ref{sigmaprev})-(\ref{bigsum}), we
obtain
\begin{equation}\label{sone}
\Sigma_1=O(\alpha_n^{-1}\log^{\rho+1}{\alpha_n^{-1}}).
\end{equation}

$\Sigma_2$ can be estimated in a similar way. We have
\begin{eqnarray} \label{stwo}
 & & \Sigma_2
 =O\left(\sum_{k=m_1+1}^\infty
 \frac{k^{\rho+1}e^{-k\alpha_n}}{1-e^{-k\alpha_n}}\right) =O\left(\sum_{k=m_1+1}^\infty
k^{\rho+1}e^{-k\alpha_n}\right) \nonumber \\
& & =O\left(\alpha_n^{-\rho-2}\int_{m_1\alpha_n}^\infty
u^{\rho+1}e^{-u}du\right)
=O(\alpha_n^{-\rho-2}(m_1\alpha_n)^{\rho+1}
e^{-m_1\alpha_n}) \nonumber \\
& & =O(\alpha_n^{-1}\log^{\rho+1}{\alpha_n^{-1}}),
\end{eqnarray}
where the last equality follows from (\ref{monen}) and
(\ref{expbeh}), while in the previous one we have used the the
asymptotic behavior of the incomplete gamma function $\Gamma(a,z)$
as $z\to\infty$ (see [1; Sect. 6.5]). The required estimate now
follows from (\ref{diff}), (\ref{sone}) and (\ref{stwo}).$\quad
\rule{2mm}{2mm}$

Our last forthcoming lemma supplies us with integral
representations for $F_b(e^{-\alpha})$ and $F_{m,b}(e^{-\alpha}),
\alpha>0$, using the Dirichlet series (\ref{dofes}) and its partial
sums
\begin{equation}\label{dmofes}
D_m(s)=\sum_{k=1}^m b_k k^{-s}, \quad s=\sigma+iy, \quad m\ge 1.
\end{equation}
The proof is based on a Mellin transform technique and can be
found in [13], [10; Lemma 2(ii)] and [2; Sect. 6.2].

\begin{lemma} For any $\alpha, \Delta>0$, we have
\begin{equation} \label{efem}
F_{m,b}(e^{-\alpha})= \frac{1}{2\pi
i}\int_{\rho+\Delta-i\infty}^{\rho+\Delta+i\infty}
\alpha^{-s}\Gamma(s)\zeta(s+1)D_m(s)ds
\end{equation}
and
\begin{equation} \label{ef}
F_b(e^{-\alpha})= \frac{1}{2\pi i}
\int_{\rho+\Delta-i\infty}^{\rho+\Delta+i\infty}
\alpha^{-s}\Gamma(s)\zeta(s+1)D(s)ds,
\end{equation}
where $D_m(s)$ and $D(s)$ are defined by (\ref{dmofes}) and
(\ref{dofes}), respectively.$\quad \rule{2mm}{2mm}$
\end{lemma}

\section{Proof of the Main Result}
 \setcounter{equation}{0}

 We apply first the Cauchy coefficient formula to (\ref{df}) using the
 circle $x=e^{-\alpha_n+i\theta}, \pi<\theta\le\pi$, as a contour
 of integration ($\alpha_n$ is determined by (\ref{alphan})). We
 obtain
 $$
 p_b(n)\mathbb{P}(X_n\le m) =\frac{e^{n\alpha_n}}{2\pi}\int_{-\pi}^\pi
f_{m,b}(e^{-\alpha_n+i\theta}) e^{-i\theta n}d\theta.
 $$
 Then, we break up the range of integration as follows:
 \begin{equation} \label{sumint}
p_b(n)\mathbb{P}(X_n\le m) =J_1(m,n) +J_2(m,n),
 \end{equation}
 where
 \begin{eqnarray} \label{jone}
 & & J_1(m,n) =\frac{e^{n\alpha_n}}{2\pi}\int_{-\delta_n}^{\delta_n}
f_{m,b}(e^{-\alpha_n+i\theta}) e^{-i\theta n}d\theta \nonumber \\
& &=\frac{e^{n\alpha_n+F_{m,b}(e^{-\alpha_n})}}{2\pi}
\int_{-\delta_n}^{\delta_n}
\frac{f_{m,b}(e^{-\alpha_n+i\theta})}{f_{m,b}(e^{-\alpha_n})}e^{-i\theta
n}d\theta,
\end{eqnarray}
\begin{equation} \label{jtwo}
J_2(m,n)
=\frac{e^{n\alpha_n+F_{m,b}(e^{-\alpha_n})}}{2\pi}\int_{\delta_n<\mid\theta\mid\le\pi}
\frac{f_{m,b}(e^{-\alpha_n+i\theta})}{f_{m,b}(e^{-\alpha_n})}e^{-i\theta
n}d\theta
\end{equation}
($\delta_n$ and $F_{m,b}(x)$ are defined by (\ref{delta}) and
(\ref{capitalfm}), respectively).

We start with an estimate for $J_1(m,n)$, expanding the integrand
of (\ref{jone}) by Taylor formula:
\begin{eqnarray}
 & & \frac{f_{m,b}(e^{-\alpha_n+i\theta})}{f_{m,b}(e^{-\alpha_n})}
 =\exp\{(e^{i\theta}-1)e^{-\alpha_n}F_{m,b}^\prime(e^{-\alpha_n})
 \nonumber \\
 & & +\frac{1}{2}(e^{i\theta}-1)^2 e^{-2\alpha_n}F_{m,b}^{\prime\prime}(e^{-\alpha_n})
 +O(\mid\theta\mid^3
 F_{m,b}^{\prime\prime\prime}(e^{-\alpha_n}))\}. \nonumber
 \end{eqnarray}
 Hence, we can rewrite (\ref{jone}) as follows:
 \begin{equation} \label{jonei}
 J_1(m,n)
 =\frac{e^{n\alpha_n}}{\sqrt{2\pi}}e^{F_{m,b}(e^{-\alpha_n})}I_n,
 \end{equation}
 where
 \begin{eqnarray}
 & & I_n=\frac{1}{\sqrt{2\pi}} \int_{-\delta_n}^{\delta_n}
 \exp\{(e^{i\theta}-1)e^{-\alpha_n}F_{m,b}^\prime(e^{-\alpha_n})
 \nonumber \\
 & & +\frac{1}{2}(e^{i\theta}-1)^2 e^{-2\alpha_n}F_{m,b}^{\prime\prime}(e^{-\alpha_n})
 +O(\mid\theta\mid^3
 F_{m,b}^{\prime\prime\prime}(e^{-\alpha_n}))-i\theta n\}d\theta.
 \nonumber
 \end{eqnarray}
 Lemma 7 shows that, for those integers $m$ satisfying (\ref{mn}), we can replace the
 derivatives $F_{m,b}^{(j)}(e^{-\alpha_n})$ by
 $F_b^{(j)}(e^{-\alpha_n}), j=1,2,3,$ at the expense of a
 negligible error term. In fact, combining Lemma 7 with
 (\ref{delta}), we have
 \begin{eqnarray}
 & & (e^{i\theta}-1)^j F_{m,b}^{(j)}(e^{-\alpha_n}) =(e^{i\theta}-1)^j
 F_b^{(j)}(e^{-\alpha_n})
 +O(\delta_n^j\alpha_n^{-j}\log^{\rho+j}{\alpha_n^{-1}}) \nonumber
 \\
 & & =(e^{i\theta}-1)^j
 F_b^{(j)}(e^{-\alpha_n}) +O(\alpha_n^{\rho
 j/3}\log^{\rho+j}{\alpha_n^{-1}}/\omega^j(n)), \quad j=1,2,
 \nonumber
 \end{eqnarray}
and
$$
O(\mid\theta\mid^3 F_{m,b}^{\prime\prime\prime}(e^{-\alpha_n})) =
O(\mid\theta\mid^3 F_b^{\prime\prime\prime}(e^{-\alpha_n}))+
 O(\alpha_n^\rho(\log^{\rho+3}{\alpha_n^{-1}})/\omega^3(n)),
$$
where the function $\omega(n)\to\infty$ as $n\to\infty$
arbitrarily slowly. Since all error terms above tend to $0$, we
obtain
\begin{eqnarray}
& & I_n =\frac{1+o(1)}{\sqrt{2\pi}} \int_{-\delta_n}^{\delta_n}
 \exp\{F_b(e^{-\alpha_n+i\theta}) -F_b(e^{-\alpha_n}) -i\theta n\}
 d\theta \nonumber \\
 & & =\frac{1+o(1)}{\sqrt{2\pi}} \int_{-\delta_n}^{\delta_n}
 \frac{f_b(e^{-\alpha_n+i\theta})}{f_b(e^{-\alpha_n})} e^{-i\theta n}d\theta, \nonumber
 \end{eqnarray}
 where the last equality follows from a similar
 Taylor expansion for $F_b(e^{-\alpha_n+i\theta})$ and the fact that
 $O(\mid\theta\mid^3 F_b^{\prime\prime\prime}(e^{-\alpha_n}))
 =O(\delta_n^3 F_b^{\prime\prime\prime}(e^{-\alpha_n}))=O(1/\omega^3(n))=o(1)$ (see (\ref{terz}) of Lemma 2
 and (\ref{delta})). Now, from
 Lemma 3 it follows that
 \begin{eqnarray}
 & & I_n \sim \frac{1}{\sqrt{2\pi}} \int_{-\delta_n}^{\delta_n}
 e^{-\theta^2 \mathcal{B}_b(e^{-\alpha_n})/2} d\theta \nonumber \\
 & & =\frac{1}{\sqrt{2\pi\mathcal{B}_b(e^{-\alpha_n})}}
 \int_{-\delta_n\sqrt{\mathcal{B}_b(e^{-\alpha_n})}}^{\delta_n\sqrt{\mathcal{B}_b(e^{-\alpha_n})}}
 e^{-y^2/2}dy \sim\frac{1}{\sqrt{2\pi\mathcal{B}_b(e^{-\alpha_n})}}
 \int_{-\infty}^\infty e^{-y^2/2}dy \nonumber \\
 & & =\frac{1}{\sqrt{\mathcal{B}_b(e^{-\alpha_n})}}, \quad
 n\to\infty. \nonumber
 \end{eqnarray}
The last asymptotic equivalence follows from (\ref{balpha}) of
Lemma 2, which implies that
$$
\delta_n\sqrt{\mathcal{B}_b(e^{-\alpha_n})} \sim
\frac{\alpha_n^{-\rho/6}}{\omega(n)} \sqrt{h(\rho+1)}\to\infty,
$$
if $\omega(n)\to\infty$ slower than $\alpha_n^{-\rho/6}$.

Substituting the asymptotic equivalence for $I_n$ into
(\ref{jonei}), we conclude that
\begin{equation} \label{joneas}
J_1(m,n)\sim
\frac{e^{n\alpha_n}}{\sqrt{2\pi\mathcal{B}(e^{-\alpha_n})}}
e^{F_{m,b}(e^{-\alpha_n})}
\end{equation}
if $m$ satisfies (\ref{mn}) as $n\to\infty$.

For the estimate of $J_2(m,n)$, we recall (\ref{jtwo}) and the
proof of (\ref{realparts}). Thus, for any real $u$, we obtain
\begin{eqnarray} \label{realpartsfm}
& & \Re{(F_{m,b}(e^{-\alpha_n+2\pi iu}))} -F_{m,b}(e^{-\alpha_n})
\le -\frac{\log{5}}{2}\sum_{k=1}^m b_k e^{-\alpha_n k}\sin^2{(\pi
uk)} \nonumber \\
& & =-\frac{\log{5}}{2}\left(\sum_{k=1}^\infty b_k e^{-\alpha_n
k}\sin^2{(\pi uk)} -\sum_{k=m+1}^\infty b_k e^{-\alpha_n
k}\sin^2{(\pi uk)}\right).
\end{eqnarray}
Consider again the sequences $m=m(n)$ and $m_1=m_1(n)$ defined by
(\ref{mn}) and (\ref{monen}), respectively. We have
 \begin{eqnarray}\label{twosums}
 & & \sum_{k=m+1}^\infty b_k e^{-\alpha_n k} \sin^2{(\pi ku)}
\le\sum_{k=m+1}^\infty b_k e^{-\alpha_n k} \nonumber \\
 & & =\sum_{k=m+1}^{m_1}b_k
e^{-\alpha_n k} +\sum_{k=m_1+1}^\infty b_k e^{-\alpha_n k}.
 \end{eqnarray}
These sums can be estimated using the argument given in the proof
of Lemma 7 (see (\ref{sigmaprev})-(\ref{stwo})). We obtain in the
same way that
\begin{equation}\label{firsts}
\sum_{k=m+1}^{m_1}b_k e^{-\alpha_n k}
=O\left(\alpha_n^\rho\sum_{k=m+1}^{m_1} b_k\right) =O((\alpha_n
m_1)^\rho) =O(\log^\rho{\alpha_n^{-1}}),
\end{equation}
\begin{equation}\label{seconds}
\sum_{k=m_1+1}^\infty b_k e^{-\alpha_n k}
=O\left(\sum_{k=m_1+1}^\infty k^\rho e^{-\alpha_n k}\right)
=O(\log^\rho{\alpha_n^{-1}})
\end{equation}
and thus, (\ref{twosums})-(\ref{seconds}) imply that
$$
\sum_{k=m+1}^\infty b_k e^{-\alpha_n k} \sin^2{(\pi ku)}
=O(\log^\rho{\alpha_n^{-1}}).
$$
Replacing the second term of the right-hand side of
(\ref{realpartsfm}) by the last $O$-estimate and applying
inequality (\ref{sn}) (see also (\ref{sinus})) to its first term,
for $\delta_n/2\pi\le\mid u\mid <1/2$, we obtain
$$
 \Re{(F_{m,b}(e^{-\alpha_n+2\pi iu}))} -F_{m,b}(e^{-\alpha_n})
\le -C_3\alpha_n^{-\epsilon_1} +O(\log^\rho{\alpha_n^{-1}}).
$$
Now, we are ready to compare the growth of (\ref{jtwo}) with that
of (\ref{joneas}) whenever $m$ satisfies (\ref{mn}). We have
\begin{eqnarray}\label{jtwoas}
& & \mid J_2(m,n)\mid \le\exp{(n\alpha_n+F_{m,b}(e^{-\alpha_n}))}
\nonumber \\
& & \times\int_{\frac{\delta_n}{2\pi}<\mid u\mid\le\frac{1}{2}}
\mid f_{m,b}(e^{-\alpha_n+2\pi iu})/f_{m,b}(e^{-\alpha_n})\mid du
\nonumber \\
& & =\exp{(n\alpha_n+ F_{m,b}(e^{-\alpha_n}))}
\int_{\frac{\delta_n}{2\pi}<\mid u\mid\le\frac{1}{2}}
(\Re{(F_{m,b}(e^{-\alpha_n+2\pi iu})} -F_{m,b}(e^{-\alpha_n})) du
\nonumber \\
& & =O(\exp{(n\alpha_n+F_{m,b}(e^{-\alpha_n})
-C_3\alpha_n^{-\epsilon_1} +O(\log^\rho{\alpha_n^{-1}}))})
\nonumber
\\
& & =O(e^{-C_3\alpha_n^{-\epsilon_1}}
\sqrt{2\pi\mathcal{B}_b(e^{-\alpha_n})}J_1(m,n)) =o(J_1(m,n)),
\end{eqnarray}
where for the last $o$-estmate we have used (\ref{balpha}). It is
now clear that (\ref{sumint}), (\ref{joneas}) and (\ref{jtwoas})
imply that
$$
p_b(n)\mathbb{P}(X_n\le m)
\sim\frac{e^{n\alpha_n}}{\sqrt{2\pi\mathcal{B}_b(e^{-\alpha_n})}}
e^{F_{m,b}(e^{-\alpha_n})}, \quad n\to\infty.
$$
Subsequent application of the asymptotic equivalence
(\ref{totalnr}) from Lemma 5 implies that
\begin{equation} \label{dfas}
\mathbb{P}(X_n\le m) \sim
\exp{\{F_{m,b}(e^{-\alpha_n})-F_b(e^{-\alpha_n})\}},
\end{equation}
where $\alpha_n$ and $m$ satisfy (\ref{alphan}) and (\ref{mn}),
respectively.

Further on we shall study the asymptotic behavior of the exponent
in (\ref{dfas}). Our analysis will be based on a generalization of
Perron formula that expresses the partial sums of a Dirichlet
series as complex integrals of the inverse Mellin type transforms
applied to the Dirichlet series itself. We shall use it in the
form given in the Supplement of [18; Sect. 3]. So, first we
represent $F_{m,b}(e^{-\alpha_n})$ using eq. (\ref{efem}) of Lemma
8 and then, we apply Perron formula to the partial sum $D_m(s)$ of
the Dirichlet series $D(s)$ (recall also (\ref{dmofes}) and
(\ref{dofes})). In this way we arrive at the following complex
integral representation: for any $\Delta>1$, we have
\begin{eqnarray} \label{doubleint}
& & F_{m,b}(e^{-\alpha_n})  \\
& & = \frac{1}{2\pi i}
\int_{\rho+\Delta-i\infty}^{\rho+\Delta+i\infty}
\alpha_n^{-s}\Gamma(s)\zeta(s+1)
\left(\frac{A(m+1)^{\rho-s}}{\rho-s} +D(s) +\Omega_m\right) ds,
\nonumber
\end{eqnarray}
where $\Omega_m=o(1), m\to\infty$. Furthermore, (\ref{dfas}) and
(\ref{doubleint}) imply that
\begin{eqnarray} \label{dfasint}
& & \mathbb{P}(X_n\le m)  \\
& & \sim \exp{\left\{-\frac{A\alpha_n^{-\rho}}{2\pi i}
\int_{\Delta-i\infty}^{\Delta+i\infty} ((m+1)\alpha_n)^{-s}
\Gamma(s+\rho) \zeta(s+\rho+1) \frac{ds}{s}\right\}}. \nonumber
\end{eqnarray}
The proofs of (\ref{doubleint}) and (\ref{dfasint}) contain some
technical details that will be given in the Appendix.

We continue with the computation of the complex integral in the
exponent of (\ref{dfasint}). The sequence $m=m(n)$ will be
specified later in a more precise way. At this moment we only
asume that it satisfies (\ref{mn}). We set in the integral of
(\ref{dfasint})
\begin{equation} \label{un}
u=u_n=(m+1)\alpha_n,
\end{equation}
and consider it as a function of $u$. First, we shall obtain its
explicit form and then we shall estimate it as $u\to\infty$ (see
(\ref{mn}) and (\ref{un})). Clearly, we can consider this integral
as the inverse Mellin transform of the function
$\Gamma(s+\rho)\zeta(s+\rho+1)/s$. For the sake of convenience, we
set
\begin{equation}\label{invmell}
H(u) =\frac{1}{2\pi i} \int_{\Delta-i\infty}^{\Delta+i\infty}
u^{-s}\Gamma(s+\rho)\zeta(s+\rho+1) \frac{ds}{s}.
\end{equation}
It is known that, for $\Re{(s)}>0$, $g_1(s)=1/s$ is the Mellin
transform of the (Heaviside-like) step function
 $$
 H_1(u)=\left\{\begin{array}{ll} 1 & \qquad  \mbox {if}\qquad 0\le u<1, \\
 0 & \qquad \mbox {if}\qquad u>1,
 \end{array}\right.
  $$
  while $g_2(s)=\Gamma(s)\zeta(s+1)$ is the
  Mellin transform of
  \begin{equation}\label{htwo}
  H_2(u) =\sum_{j=1}^\infty\frac{e^{-ju}}{j} =-\log{(1-e^{-u})}
  \end{equation}
  (see e.g. [7; Appendix B.7]). Next, for $\Delta>1$, we apply
  formula (6.1.14) from [5] with $\alpha=0$ and $\beta=\rho-1$. We
  obtain
  \begin{eqnarray} \label{hu}
  & & H(u) =u^\alpha\int_0^\infty y^\beta H_1(u/y)H_2(y)dy
  =-\int_u^\infty y^{\rho-1} \log{(1-e^{-y})} dy \nonumber \\
  & & = \int_u^\infty y^{\rho-1} e^{-y} dy +R(u)=\Gamma(\rho,u) +R(u),
  \end{eqnarray}
  where $\Gamma(\rho,u)$ denotes the incomplete gamma function,
  while $R(u)$ is the error term given by
  $$
  R(u)=\int_u^\infty
  y^{\rho-1} \left(\sum_{j=2}^\infty\frac{e^{-jy}}{j}\right)dy.
  $$
  It is easily estimated as follows:
  $$
  R(u) \le\frac{1}{2(1-e^{-u})}
  \int_u^\infty
  y^{\rho-1} e^{-2y}dy=O(e^{-u}\Gamma(\rho,u)), \quad u\to\infty.
  $$
  Combining this estimate with (\ref{dfasint})-(\ref{invmell}) and (\ref{hu}) and applying
  the asymptotic $\Gamma(\rho,u_n)=u_n^{\rho-1}e^{-u_n}(1+O(1/u_n))$ of the incomplete gamma function (see again [1;
  Sect. 6.5]), we obtain
  \begin{eqnarray} \label{dfasfin}
  & & \mathbb{P}(X_n\le m) \sim
  \exp{\{-A\alpha_n^{-\rho}(u_n^{\rho-1}e^{-u_n}(1+O(1/u_n))
  +O(u_n^{\rho-1}e^{-2u_n}))\}} \nonumber \\
  & & = \exp{\{-A\alpha_n^{-1} m^{\rho-1}e^{-m\alpha_n}
  (1+O(1/m\alpha_n))\}} \nonumber \\
  & & =\exp{\{-A\alpha_n^{-1} m^{\rho-1}e^{-m\alpha_n}
  (1+O(1/\log{\alpha_n^{-1}}))\}},
  \end{eqnarray}
  where for the last equality we have used again (\ref{mn}).
  It is now clear that $\mathbb{P}(X_n\le m)$ converges to the
  distribution function $e^{-e^{-t}}, -\infty<t<\infty$, if
  $m=m(n)$ satisfies
  $$
  -m\alpha_n +(\rho-1)\log{m} +\log{(A\alpha_n^{-1})} =-t+o(1)
  $$
  as $n\to\infty$. From this we deduce
  \begin{equation} \label{emlogm}
  m=\alpha_n^{-1}\log{\alpha_n^{-1}} +(\rho-1)\alpha_n^{-1}\log{m}
  +(\log{A}+t)\alpha_n^{-1} +o(\alpha_n^{-1}),
  \end{equation}
  which in turn implies that
  \begin{eqnarray}
  & & \log{m} =\log{(\alpha_n^{-1}\log{\alpha_n^{-1}})} \nonumber \\
  & & +\log{\left(1+\frac{\log{A}+t}{\log{\alpha_n^{-1}}}
  +(\rho-1)\frac{\log{m}}{\log{\alpha_n^{-1}}}+o(1/\log{\alpha_n^{-1}})\right)} \nonumber \\
  & & = \log{\alpha_n^{-1}} +\log{\log{\alpha_n^{-1}}} \nonumber
  \\
  & & +\log{\left(1+\frac{\log{A}+t}{\log{\alpha_n^{-1}}}
  +(\rho-1)\frac{\log{\alpha_n^{-1}}+\log{\log{\alpha_n^{-1}}}+O(1)}
  {\log{\alpha_n^{-1}}}\right)} \nonumber \\
  & & =\log{\alpha_n^{-1}} +\log{\log{\alpha_n^{-1}}}
  +\log{\left(1+(\rho-1)+
  O\left(\frac{\log{\log{\alpha_n^{-1}}}}{\log{\alpha_n^{-1}}}\right)\right)}
  \nonumber \\
  & & =\log{\alpha_n^{-1}} +\log{\log{\alpha_n^{-1}}} +\log{\rho}
  +O\left(\frac{\log{\log{\alpha_n^{-1}}}}{\log{\alpha_n^{-1}}}\right).
  \nonumber
  \end{eqnarray}
 Hence, (\ref{emlogm}) becomes
  \begin{eqnarray}
  & & m=\rho\alpha_n^{-1}\log{\alpha_n^{-1}}
  +(\rho-1)\alpha_n^{-1}\log{\log{\alpha_n^{-1}}} \nonumber \\
  & & +\alpha_n^{-1}(\rho-1)\log{\rho} +(\log{A}+t)\alpha_n^{-1}
  +O\left(\alpha_n^{-1}\frac{\log{\log{\alpha_n^{-1}}}}{\log{\alpha_n^{-1}}}\right).
  \nonumber
  \end{eqnarray}
  Replacing now this value of $m$ into (\ref{dfasfin}) and using the continuity of the
  distribution function $e^{-e^{-t}}, -\infty<t<\infty$, we obtain
  \begin{eqnarray} \label{lthalpha}
  & & \mathbb{P}(X_n\le m) \nonumber \\
  & & =\mathbb{P}(\alpha_n X_n -\rho\log{\alpha_n^{-1}} -(\rho-1)\log{\log{\alpha_n^{-1}}}
  -(\rho-1)\log{\rho} -\log{A}+o(1)\le t) \nonumber \\
  & & \to e^{-e^{-t}}, \quad n\to\infty.
  \end{eqnarray}
 To complete the proof of the theorem, it
 remains to justify the normalization for $X_n$ stated in (\ref{limth}). We have to show that the
 sequence
 $\alpha_n, n\ge 1$, in (\ref{lthalpha}) can be replaced by $a(n)=a(n;\rho,A), n\ge
 1$ (see (\ref{norm})). So, we recall first (\ref{alphan}) and
 notice that taking logarithms from its both sides, we
 easily obtain
 \begin{eqnarray}\label{logan}
 & & \log{\alpha_n^{-1}} =\mid\log{a(n)}\mid +O(n^{-\frac{\rho}{\rho+1}}),
 \nonumber \\
 & & \log{\log{\alpha_n^{-1}}} =\log{\mid\log{a(n)\mid}} +O(n^{-\frac{\rho}{\rho+1}}/\log{n}).
 \end{eqnarray}
 Next, we set $Z_n:=\alpha_n X_n+z_n, Y_n:=a(n) X_n+y_n$, where
 $z_n:=\rho\log{\alpha_n}-(\rho-1)\log{\log{\alpha_n^{-1}}}
 -(\rho-1)\log{\rho}-\log{A}, y_n=\rho\log{a(n)}
 -(\rho-1)\log{\mid\log{a(n)}\mid}-(\rho-1)\log{\rho}-\log{A}, n\ge 1$. Furthermore, (\ref{lthalpha}) can be written
 in a shorter way as follows:
 \begin{equation} \label{conv}
  G_{n,Z}(t):=\mathbb{P}(Z_n\le t)\to e^{-e^{-t}}, \quad n\to\infty.
 \end{equation}
 So, we have to prove the same convergence for
 $G_{n,Y}(t):=\mathbb{P}(Y_n\le t), n\ge 1$. It is easy to verify that the
 above representations for $Y_n$ and $Z_n$ imply that
 \begin{equation}\label{repres}
 Y_n=(1+\eta_n)Z_n+e_n,
 \end{equation}
 where $\eta_n:=a(n)/\alpha_n-1, e_n:=y_n-(a(n)/\alpha_n)z_n, n\ge
 1$. From (\ref{norm}), (\ref{alphan}) and (\ref{logan}) we obtain
 the estimates $\eta_n=O(n^{-\frac{\rho}{\rho+1}})$ and
 $e_n=O(n^{-\frac{\rho}{\rho+1}}\log{n})$ as $n\to\infty$. Moreover, (\ref{repres}) implies that
 \begin{eqnarray}
 & & G_{n,Y}(t) =\mathbb{P}\left(Z_n\le \frac{t-e_n}{1+\eta_n}\right)
 =\mathbb{P}\left(Z_n\le t-\frac{t\eta_n+e_n}{1+\eta_n}\right) \nonumber \\
 & & =G_{n,Z}\left(t-\frac{t\eta_n+e_n}{1+\eta_n}\right), \quad
 n\ge 1. \nonumber
 \end{eqnarray}
 Taking now an arbitrary $\eta>0$ and $n$ enough large so that
 $-\eta<(t\eta_n+e_n)/(1+\eta_n)<\eta$, for fixed $t$, we obtain
 $$
 G_{n,Z}(t-\eta)\le G_{n,Y}(t)\le G_{n,Z}(t+\eta).
 $$
 Letting $n\to\infty$ in the above inequalities, from (\ref{conv})
 we find that
 $$
 e^{-e^{-(t-\eta)}}\le\lim\inf_{n\to\infty}G_{n,Y}(t)
 \le\lim\sup_{n\to\infty}G_{n,Y}(t) \le e^{-e^{-(t+\eta)}}
 $$
 for all $\eta>0$. Letting now $\eta\to 0^+$, the required result stated in (\ref{limth})
 follows from the continuity of the distribution function
 $e^{-e^{-t}}$. $\quad \rule{2mm}{2mm}$

\section*{Appendix}
\renewcommand{\theequation}{A.\arabic{equation}}
\setcounter{equation}{0}

{\it Proof of (\ref{doubleint}).} First, we recall eq.
(\ref{efem}) of Lemma 8. Our goal is to represent the $m$th
partial sum $D_m(s)$, defined by (\ref{dmofes}), using the
inversion formula given by Thm. 3.1 in the Supplement of [18] (see
also formula (3.4) there). Instead of $D(s)$, we shall consider
now the Dirichelet series $\tilde{D}(s)=D(s+\rho-1)$ (see
(\ref{dtil}) and condition $(M_1)$). It converges absolutely for
$\Re{(s)}=\sigma>1$. Lemma 6 implies that the coefficients of
$D(s+\rho-1)$ satisfy
$$
b_k k^{-\rho+1} =o(k^{\rho})k^{-\rho+1} =o(k)<\tilde{c}k,
$$
for some constant $\tilde{c}>0$ and all $k$ (in other words,
$\Phi(x)=x$ in Thm. 3.1 of [18; Supplement]). Furthermore, from
condition $(M_1)$ it follows that
$$
\sum_{k=1}^\infty b_k k^{-\rho+1} k^{-\sigma} =\frac{A}{\sigma-1}
+\phi(s),
$$
where $\phi(s)$ denotes a function which is analytic for
$\sigma\ge -C_0$. Hence
$$
\sum_{k=1}^\infty b_k k^{-\rho+1} k^{-\sigma} =
O((\sigma-1)^{-1}), \quad \sigma\to 1^+.
$$
So, the conditions of Thm. 3.1 (the Supplement of [18]) are
satisfied and by its second part we conclude that, for large
enough $T>0$, $\Delta>1$ and $d>0$, we have
\begin{eqnarray} \label{dmint}
& & D_m(w+\rho-1) +\frac{1}{2}b_{m+1}(m+1)^{-\rho+1}(m+1)^{-w}
\nonumber \\
& & =\frac{1}{2\pi i} \int_{d-iT}^{d+iT} D(w+z+\rho-1)
\frac{(m+1)^z}{z} dz
  +O\left(\frac{m^d}{(\Delta+d)T}\right) \nonumber \\
 & & + O\left(\frac{m^{1-\Delta}\log{m}}{T}\right), \quad w=1+\Delta+iy, -\infty<y<\infty.
\end{eqnarray}
Lemma 6 implies that the second term in the left-hand side of
(\ref{dmint}) is $o(m^{-\Delta})$ as $m\to\infty$. To compute the
integral in the right-hand side of (\ref{dmint}), we set
$d=1/\log{m}, m\ge 2$, and use a contour integral around the
rectangle $d-iT, d+iT, -C_0-\rho-\Delta+iT, -C_0-\rho-\Delta-iT$.
Using condition $(M_2)$, we estimate the integral over the end
segment $(-C_0-\rho-\Delta+iT,-C_0-\rho-\Delta-iT)$ by
$O(T^{C_1+1} m^{-C_0-\rho-\Delta})$. Hence, it tends to $0$ as
$m,T\to\infty$, provided $T=o(m^{(C_0+\rho+\Delta)/(C_1+1)})$. The
integrals on the segments $(-C_0-\rho-\Delta+iT,d+iT)$ and
$(-C_0-\rho-\Delta-iT,d-iT)$ are easily estimated. By condition
$(M_2)$ and the choice of $d$ both are of order
$$
O\left(T^{C_1-1}\int_{-C_0-\rho-\Delta}^{1/\log{m}}(m+1)^\sigma
d\sigma\right)
=O\left(\frac{T^{C_1-1}m^{1/\log{m}}}{\log{m}}\right)
=O\left(\frac{T^{C_1-1}}{\log{m}}\right).
$$
 Therefore, we
conclude that $T$ should satisfy
\begin{equation}\label{mttwo}
 T =\left\{\begin{array}{ll} o(m^{(C_0+\rho+\Delta)/(C_1+1)}) & \qquad  \mbox {if}\qquad C_1\le 1, \\
 o((\log{m})^{1/(C_1-1)}) & \qquad \mbox {if}\qquad C_1>1.
 \end{array}\right.
 \end{equation}
Further on we shall assume that $T=T(m)\to\infty$ as $m\to\infty$
and $m$ and $T$ satisfy (\ref{mttwo}), where the constants $C_0$,
$\rho$ and $C_1$ are defined by conditions $(M_1)$ and $(M_2)$ and
$\Delta>1$ is fixed. Thus, all integrals on the end segments
except the integral on $(d-iT,d+iT)$ are close to $0$ for enough
large $m$ and $T$. The same obviously holds for both $O$-estimates
in the right-hand side of (\ref{dmint}). So, we can compute
$D_m(w+\rho-1)$ summing up the residues of the integrand in
(\ref{dmint}). Inside the contour of integration it has only two
simple poles: at $z=1-w$ and $z=0$. Thus, we obtain
\begin{eqnarray}\label{dmlast}
& & D_m(w+\rho-1) =\frac{A(m+1)^{1-w}}{1-w} +D(w+\rho-1)
+\Omega_m, \\
& & w=1+\Delta+iy, -\infty<y<\infty, \nonumber
\end{eqnarray}
 where $\Omega_m$ equals the sum of all negligible
 terms described above. Clearly,
 \begin{equation}\label{bigomega}
\Omega_m\to 0, \quad m\to\infty,
\end{equation}
for $m$ and $T$ satisfying (\ref{mttwo}). Setting in
(\ref{dmlast}) $w=s-\rho+1$ and substituting this expression into
(\ref{efem}), we arrive at (\ref{doubleint}). $\quad
\rule{2mm}{2mm}$

{\it Proof of (\ref{dfasint}).} (\ref{doubleint}) and (\ref{ef})
imply that
\begin{eqnarray}
& & F_{m,b}(e^{-\alpha_n}) =\frac{A}{2\pi i}
\int_{\rho+\Delta-i\infty}^{\rho+\Delta+i\infty} \alpha_n^{-s}
\Gamma(s)\zeta(s+1) \frac{(m+1)^{\rho-s}}{\rho-s}ds \nonumber \\
& & +F_b(e^{-\alpha_n}) +\frac{A\Omega_m}{2\pi i}
\int_{\rho+\Delta-i\infty}^{\rho+\Delta+i\infty} \alpha_n^{-s}
\Gamma(s)\zeta(s+1)ds.  \nonumber
\end{eqnarray}
The last integral represents an inverse Mellin transform whose
original (see (\ref{htwo}) and [7; Appendix B.7]) is
$H_2(\alpha_n)=-\log{(1-e^{-\alpha_n})}=O(-\log{\alpha_n})$, as
$n\to\infty$. Hence assumption (\ref{mn}) and (\ref{bigomega})
imply that
\begin{eqnarray}\label{fmbint}
& & F_{m,b}(e^{-\alpha_n}) =\frac{A}{2\pi i}
\int_{\rho+\Delta-i\infty}^{\rho+\Delta+i\infty} \alpha_n^{-s}
\Gamma(s)\zeta(s+1) \frac{(m+1)^{\rho-s}}{\rho-s}ds \nonumber \\
& & +F_b(e^{-\alpha_n}) +o(-\log{\alpha_n}). \nonumber
\end{eqnarray}
To obtain (\ref{dfasint}) it is enough to replace this expression
into (\ref{dfas}), change the variable $s$ in the above integral
by $s+\rho$ and observe that
$-\log{\alpha_n}=o(\alpha_n^{-\rho})$. $\quad \rule{2mm}{2mm}$

\section*{Acknowledgements}

 The author would like to thank the referee for his valuable comments
 and especially for indicating some defects contained in previous versions of the paper.

\end{document}